\input gtmacros
\input amsnames
\input amstex
%
\catcode`\@=12        
\input gtmonout
\volumenumber{1}
\volumeyear{1998}
\volumename{The Epstein birthday schrift}
\pagenumbers{99}{116}
\papernumber{4}
\received{16 November 1997}
\published{21 October 1998}

%
\let\\\par
\def\topmatter{\relax}
\def\endtopmatter{\maketitlepage}
\let\gttitle\title
\def\title#1\endtitle{\gttitle{#1}}
\let\gtauthor\author
\def\author#1\endauthor{\gtauthor{#1}}
\let\gtaddress\address
\def\address#1\endaddress{\gtaddress{#1}}
\let\gtemail\email
\def\email#1\endemail{\gtemail{#1}}
\def\subjclass#1\endsubjclass{\primaryclass{#1}}
\let\gtkeywords\keywords
\def\keywords#1\endkeywords{\gtkeywords{#1}}
\def\heading#1\endheading{{\def\S##1{\relax}\def\\{\relax\ignorespaces}
    \section{#1}}}
\def\head#1\endhead{\heading#1\endheading}

\def\subhead#1\endsubhead{\sh{#1}}
\def\subsubhead#1\endsubsubhead{\sh{#1}}
\def\specialhead#1\endspecialhead{\sh{#1}}
\def\demo#1{\rk{#1}\ignorespaces}
\def\enddemo{\ppar}
\let\remark\demo
\def\endremark{}
\let\definition\demo
\def\enddefinition{\ppar}
\let\example\demo
\def\endexample{\ppar}
\def\qed{\ifmmode\quad\sq\else\hbox{}\hfill$\sq$\par\medskip\goodbreak\rm\fi}  
\def\proclaim#1{\rk{#1}\sl\ignorespaces}
\def\endproclaim{\rm\ppar}
\def\cite#1{[#1]}
\newcount\itemnumber
\def\roster{\items\itemnumber=1}
\def\endroster{\enditems}
\let\itemold\item
\def\item{\itemold{{\rm(\number\itemnumber)}}%
\global\advance\itemnumber by 1\ignorespaces}
\def\S{section~\ignorespaces}  
\def\date#1\enddate{\relax}
\def\thanks#1\endthanks{\relax}   
\def\dedicatory#1\enddedicatory{\relax}  
\let\footnote\plainfootnote
\def\Refs{\ppar{\large\bf References}\ppar\bgroup\leftskip=25pt
\frenchspacing\parskip=3pt plus2pt\small}       
\def\endRefs{\egroup}
\def\widestnumber#1#2{\relax}
\def\endrefitem{}
\def\refdef#1#2#3{\def#1{\endrefitem#2\def\endrefitem{#3}}}
\def\ref{\par}
\def\endref{\endrefitem\par\def\endrefitem{}}
\refdef\key{\noindent\llap\bgroup[}{]\stdspace\egroup}
\refdef\no{\noindent\llap\bgroup[}{]\stdspace\egroup}
\refdef\by{\bf}{\rm, }
\refdef\manyby{\bf}{\rm, }
\refdef\paper{\it}{\rm, }
\refdef\book{\it}{\rm, }
\refdef\jour{}{ }
\refdef\vol{}{ }
\refdef\yr{$(}{)$ }
\refdef\ed{(}{, editor) }
\refdef\publ{}{ }
\refdef\inbook{from: ``}{'', }
\refdef\pages{}{ }
\refdef\page{}{ }
\refdef\paperinfo{}{ }
\refdef\bookinfo{}{ }
\refdef\publaddr{}{ }
\refdef\eds{(}{, editors) }
\refdef\bysame{\hbox to 3 em{\hrulefill}\thinspace,}{ }
\refdef\toappear{(to appear)}{ }
\refdef\issue{no.\ }{ }
\newcount\refnumber\refnumber=1
\def\refkey#1{\expandafter\xdef\csname cite#1\endcsname{\number\refnumber}%
\global\advance\refnumber by 1}
\def\cite#1{[\csname cite#1\endcsname]}
\def\Cite#1{\csname cite#1\endcsname}  
\def\key#1{\noindent\llap{[\csname cite#1\endcsname]\ \ }}

\refkey {Al} 
\refkey {B}
\refkey {BH}
\refkey {BHi}
\refkey {BM}
\refkey {BS}
\refkey {Ge}
\refkey {Gr}
\refkey {Ha}
\refkey {Hi}
\refkey {HNN}
\refkey {HW}
\refkey {LS}
\refkey {M}
\refkey {Me}
\refkey {Ol}
\refkey {S}
\refkey {Se}
\refkey {W1}
\refkey {W2}
\refkey {W3}

\def\G{\Gamma}
\def\c{CAT$(0)$ }
\def\e{\varepsilon}
\def\a{\alpha}
\def\b{\beta}
\def\-{\overline}
\def\A{\Cal A}

\def\g{\gamma}

\def\onto{\twoheadrightarrow}
\def\gp{\text{\rm{gp}}}
\def\<{\langle}
\def\>{\rangle}


\topmatter
\title{Controlled embeddings into groups that have\\
no non-trivial finite quotients}
\endtitle 
\shorttitle{Controlled embeddings into groups}

\author{Martin R Bridson}\endauthor

\thanks{
This work was supported by an EPSRC
Advanced Fellowship and NSF grant 9401362.} 
\endthanks

\keywords
Finite quotients, embeddings,
non-positive curvature
\endkeywords

\address{Mathematical Institute, 24--29 St Giles',
 Oxford, OX1 3LB}\endaddress

\email{bridson@maths.ox.ac.uk}\endemail

\subjclass
20E26, 20E06, 53C70;
20F32, 20F06
\endsubjclass

\abstract{If a class of
finitely generated
groups $\Cal G$ is closed under isometric amalgamations
along free subgroups, then  every $G\in \Cal G$
can be quasi-isometrically embedded in a group $\widehat G\in\Cal G$
that has no proper subgroups of finite index. 

Every compact, connected,
non-positively curved space $X$ admits
an isometric embedding
into a compact, connected,
non-positively curved space $\overline X$
such that $\overline X$ has no non-trivial finite-sheeted coverings.}
\endabstract

\asciiabstract{If a class of
finitely generated
groups Curly(G) is closed under isometric amalgamations
along free subgroups, then  every G in Curly(G)
can be quasi-isometrically embedded in a group Hat(G) in Curly(G)
that has no proper subgroups of finite index. 

Every compact, connected,
non-positively curved space X admits
an isometric embedding
into a compact, connected,
non-positively curved space Overline(X)
such that Overline(X) has no non-trivial finite-sheeted coverings.}

\endtopmatter
\document

David Epstein's lucid writings, particularly those
on automatic groups, had a strong influence on
me when I was a graduate student. Since then,
during many  hours of enjoyable conversation, I have
continued to benefit from his great insight into mathematics.
It was therefore a
great pleasure to speak at his birthday
celebration and it is an equal pleasure to
write an article for this volume. 
\sectionnumber=-1
\heading{Introduction}\endheading
In this
article  I shall address  the following general question:
given a finitely generated
group $G$ that satisfies certain desirable properties,
when can one embed $G$ into a group which retains these
desirable properties but does not have any 
non-trivial finite
quotients? My interest in this question arises from
a geometric problem that is the subject of Theorem~C.

Our discussion begins with a general embedding theorem
which is  similar to results that 
were proved in the wake of
the landmark paper by Higman, Neumann and Neumann \cite{HNN}.
The novel element in   the result presented here  is that
we control  the geometry of  the embedding. 

\proclaim{Theorem A} Let $\Cal G$ be a class of 
finitely generated groups.
If  $\Cal G$  is closed under  the operation of
isometric amalgamation 
along finitely
generated free groups,  
then   every   $G\in\Cal G$
can be quasi-isometrically
embedded in a group $\widehat G\in\Cal G$
that has no proper subgroups of finite index.
\endproclaim

The definition of  isometric amalgamation is given in
Section 1.
There are various interesting classes
of groups that are closed under 
amalgamations along arbitrary finitely generated free
groups, for example 
the class of all finitely presented groups, 
 groups
of type $F_n$, and groups of a given (cohomological or geometric)
dimension $n\ge 2$. The benefit of 
restricting the geometry of the amalgamation
becomes apparent when the defining
properties of $\Cal G$ are more geometric in nature.
For example, the class of groups which satisfy a 
polynomial isoperimetric inequality
 is   not
closed under the   operation of
 amalgamation 
along arbitrary finitely
generated free groups 
(or indeed along
quasi-isometrically embedded free groups), but it
is closed under  
 amalgamation 
along  isometrically embedded subgroups (Corollary 4.2). 

A refinement 
of the proof of Theorem A  yields:

\proclaim{Theorem B} Every finitely presented
group $G$ can be embedded in a finitely
presented group $\widehat G$
that has no non-trivial finite quotients  and whose Dehn
function $f_{\widehat G}$   satisfies:
$$
f_{\widehat G}(n)\le nf_G(n).
$$
One can (simultaneously)  arrange for the isodiametric function of
$\widehat G$ to be no greater than that of $G$.
\endproclaim

Theorem A does not apply directly to the class of groups
that arise as
fundamental groups of compact non-positively curved
spaces. \fnote{Throughout this
article we use the term `non-positive curvature'
in the sense of A.D.~Alexandrov  \cite{BH}.}
 Nevertheless, using a more subtle argument  based
on the same blueprint of proof, in Section 3 we shall prove
the following theorem. 
(We say that a covering $\widehat Z\to Z$ is `non-trivial' 
if $\widehat Z$ is connected and
$\widehat Z\to Z$ is not a homeomorphism.)

\proclaim{Theorem C} Every compact, connected,
non-positively curved space $X$ admits
an isometric embedding
into a compact, connected,
non-positively curved space $\overline X$
such that $\overline X$ has no non-trivial finite-sheeted coverings.
If $X$ is a polyhedral complex of dimension $n\ge 2$, then
one can arrange for $\overline X$ to be a complex of the same 
dimension.
\endproclaim

Any
local isometry between compact non-positively curved
spaces induces an injection on fundamental groups 
[\Cite{BH}, II.4], so in the notation of Theorem C we have
$\pi_1X\hookrightarrow\pi_1\overline X$. Since $\overline X$
has no non-trivial  finite-sheeted coverings, $\pi_1\overline X$
has no proper subgroups of finite index.
Thus
Theorem C gives a solution to our general embedding problem
for  the class of groups
that arise as
fundamental groups of compact non-positively curved
spaces.  An extension of Theorem C yields the
corresponding result for groups that act properly and
cocompactly on CAT$(0)$ spaces (3.6).

The fundamental groups of the most classical
examples of non-positively curved spaces, quotients of
symmetric spaces of non-compact type,
are residually finite. In 1995
Dani Wise
produced the
first examples of  compact non-positively
curved spaces whose fundamental groups have no non-trivial
finite quotients \cite{W3}.  He
also constructed semihyperbolic
groups that are not virtually torison free, cf~(3.7). 
Subsequently, Burger and Mozes \cite{BM}  constructed
 compact non-positively
curved 2-complexes whose fundamental groups are simple.   
Fundamental groups of compact negatively curved
spaces, on the other hand, are never simple \cite{Gr}, \cite{Ol}. 

One might hope to prove an 
analogue of Theorem A in which the enveloping
group $\widehat G$ is simple. However the techniques
described in this article
are clearly inadequate in this regard. Indeed,
  finitely presented simple groups have solvable
word problems   and hence so do  their finitely presented
subgroups. Thus if one wishes to embed a given finitely
presented group $G$ into a finitely presented 
simple group, then one must
make essential use of the fact that $G$ has a solvable
word problem. Higman conjectures that the solvability
of the word problem is the only obstruction to the
existence of such an embedding \cite{Hi} (cf~\cite{BHi},
\cite{S}).

This article is organized as follows. In Section 1 we describe
some examples of groups that are not residually finite and
define isometric amalgamation. In Section 2 we prove Theorem A.
In Section 3 we discuss spaces of non-positive curvature
and prove Theorem C. In Section 4 we examine the effect of
isometric amalgamations on isoperimetric and isodiametric
inequalities and prove Theorem B.

This article grew out of a lecture which I gave at the
conference on Geometric Group Theory at Canberra in
July 1996.  I would like to thank the organizers of that
conference. I would particularly
like to thank Chuck Miller for arranging my visit and
for welcoming me so warmly.

\eject
\heading{Residual finiteness and isometric amalgamation}
\endheading
 
A group $G$ is said to be {\it residually finite}
if for every non-trivial element $g\in G$ there is a 
finite group $Q$ and an epimorphism  $\phi\co G\onto Q$
such that $\phi(g)\neq 1$.
As a first step towards producing groups with no finite
quotients, we must gather a supply of groups that
are not residually finite. 
The Hopf property
provides a useful tool in this regard.
A group   $H$ is said to be
{\it Hopfian} if every epimorphism $H\onto H$ is an
isomorphism --- in other words, if $N\subset H$ is
normal and $H/N\cong H$ then $N=\{1\}$. 

The following result was first proved by Malcev \cite{M}.

\proclaim{1.1 Proposition} If  a finitely
generated group is residually finite
then it is Hopfian.
\endproclaim

\demo{Proof} Let $G$ be a finitely generated
group and suppose that there is an epimorphism
$\phi \co  G\to G$ with non-trivial kernel. We fix
$g_0\in \ker\phi \smallsetminus \{1\}$ and for every
  $n>0$ we choose $g_n\in G$
such that $\phi^n(g_n) = g_0$.

If there were a finite
group $Q$ and a homomorphism $p\co G\to Q$ such that
$p(g_0)\ne 1$, then all of the maps $\phi_n := p\phi^n$
would be distinct, because
 $\phi_n(g_n)\ne 1$ whereas $\phi_m(g_n)=1$
if $m>n$. But there are
only finitely many homomorphisms from any finitely
generated group to any finite group (because the images
of the generators determine the map).
\qed

\example{1.2 Examples} The following group  was discovered
by Baumslag and Solitar \cite{BS}:
$$
\text{\rm{BS}}(2,3)=\langle a, t \mid t^{-1}a^2t = a^3\rangle.
$$
The map $a\mapsto a^2, t \mapsto t$ is onto: $a$ is
in the image  because
$a=a^3 a^{-2} = (t^{-1}a^2t)a^{-2}$. However this map
is not an isomorphism: $[a, t^{-1}at]$ is a non-trivial
element of the kernel. 
Meier \cite{Me} noticed that the salient features of
this example are present in many 
other HNN extensions of abelian groups. Some
of these groups  were later studied by Wise \cite{W1}, among
them    
$$
T(n)= \langle a, b, t_a, t_b \mid [a,b]=1,\,
t_a^{-1}at_a= (ab)^n,\, t_b^{-1}bt_b=(ab)^n \rangle,
$$ 
which is the fundamental group of a compact non-positively
curved 2-complex (see (3.1)). If $n\ge 2$ then
certain non-trivial commutators, for example
$g_0=[t_a (ab) t_a^{-1},b]$,  
lie in the
kernel of the epimorphism $T(n)\onto T(n)$
given by $a\mapsto a^n, b\mapsto b^n, t_a\mapsto t_a,
t_b \mapsto t_b$.  
The proof of (1.1) shows that
$g_0$ has trivial image in every finite
quotient of $T(n)$.
\endexample

\definition{1.3 Definition of Isometric Amalgamation}
Let $H\subset G$ be a pair of 
groups with fixed finite generating
sets.  If, in   the corresponding
word metrics, $d_G(h,h')=d_H(h,h')$ for all $h,h'\in H$, 
then we say that $H$ is  {\it isometrically embedded} in $G$.

Consider a finite graph of groups (in the sense of
Serre \cite{Se}). If one can
choose  finite generating sets for the
vertex groups $G_i$ and the edge groups $H_{i,j}$ 
such that  
the inclusions of the edge groups are all
isometric embeddings, then we say that the fundamental group
$\G$
of the graph of groups is obtained by an {\it isometric
amalgamation of the $G_i$ along the $H_{i,j}$} or, more
briefly, {\it $\G$ is an isometric amalgam of the $G_i$}.

Note that, with respect to the natural
choice of generators, all of the vertex and edge groups
are isometrically embedded in the amalgam.
Note also that,
even in the basic  cases of HNN extensions
and amalgamated free products, the above definition
is more stringent than simply requiring that  for each
$i,j$
there exist  choices
of generators (depending on $i,j$) with respect
to which
$H_{i,j}\hookrightarrow G_i$ is an isometric embedding. 

Free products of finitely generated groups are
(trivial) examples of isometric amalgams. One can
also obtain both $G\times\Bbb Z$ and $G\ast\Bbb Z$
from $G$ by   isometric
amalgamations: each is the fundamental group of
a graph of groups with one vertex group $G$ and one edge 
group; to obtain $G\times\Bbb Z$ one takes $G$ as
edge group and uses
the identity map  as the inclusions; to obtain $G\ast\Bbb Z$
one takes the edge group to be trivial.
\enddefinition

\proclaim{1.4 Lemma} Let $\Cal G$ be as in Theorem A and let
$T(n)$ be as in (1.2). If
$G\in\Cal G$ then
 $G\ast T(n)\in\Cal G$.
\endproclaim

\demo{Proof} Fix a finite generating set $\Cal S$
for $G$. As above $G\ast\Bbb Z\in\Cal G$; let
$a$ be a generator of the $\Bbb Z$ free factor.
The cyclic subgroup generated by $a$ is isometrically
embedded with respect to the generating
system $\Cal S\cup\{a\}$.
 We add a further stable letter $b$ that commutes
with $a$, thus obtaining $G\ast\Bbb Z^2\in\Cal G$.

With respect to $\Cal S\cup\{a,b,(ab)^n\}$, the
cyclic subgroups generated by $a,b$ and $(ab)^n$
are all isometrically embedded. 
Thus $G\ast T(n)$ can be obtained from $G\ast\Bbb Z^2$
by an isometric amalgamation: the underlying
graph of  groups has one vertex
group, $G\ast\Bbb Z^2$,  there are two edges
in the graph and both edge groups are cyclic; the 
homomorphism at one end of each edge sends the generator to
 $ (ab)^n$, and  the maps at  the other ends are onto
 $\langle a\>$ and $\< b\>$ respectively. 
\qed

\heading{The proof of Theorem A}\endheading

In order to clarify the exposition, 
we shall first  prove a simplified version of Theorem A in
which we do not examine the geometry of
the amalgamations involved.

\proclaim{2.1 Lemma} 
Let $\Cal G$ be a class of groups that 
 is closed under  the operation of amalgamation along finitely
generated free groups. If $G\in\Cal G$ 
is finitely generated, then it can
be embedded in a finitely generated group 
 $\widehat G\in\Cal G$
that has no proper subgroups of finite index.
\endproclaim

\demo{Proof} The following proof is chosen with Theorem A in
mind (shorter proofs exist). A similar construction was
used in \cite{W3}.

{\rk{Step 0}} Replacing $G$ by $G_0=G\ast T(n)$ if necessary, we
may assume that   $G$ contains an element of infinite order
 $g_0\in G$   whose image in every
finite quotient of $G_0$ is trivial (see (1.2)). Let $\{b_1,\dots, b_n\}$
be a generating set for $G_0$. We replace
$G_0$ by   $G_1= G_0\ast\Bbb Z$, and take as
generators $\Cal A':=\{t, b_1t,\dots,b_nt\}$,
where $t$ generates the free factor $\Bbb Z$. We relabel
the generators $\Cal A'=\{a_0,\dots,a_n\}$.

{\rk{Step 1}} We take an HNN
extension of $G_1$ with $n$ stable letters: 
$$
E_1=\langle G_1, s_0, \dots, s_n \mid s_i^{-1}a_is_i = 
g_0^{p_i}, \ i=0,\dots,n\rangle.
$$
where the $p_i$ are any non-zero integers.
Now, since each $a_i$ is conjugate to a power of $g_0$
in $E_1$, the only generators of $E_1$ that can survive
in any finite quotient are the $s_i$. 
However, since there is an obvious
retraction of $E_1$ onto the free subgroup generated
by the $s_i$, the group $E_1$ still has plenty of finite
quotients.  

{\rk{Step 2}}
We  repeat the extension process, this time introducing
stable letters $\tau_i$ to make the generators $s_i$
conjugate to $g_0$:
$$
E_2=\langle E_1, \tau_0, \dots, \tau_n \mid \tau_i^{-1}s_i\tau_i 
= g_0, \ i=0,\dots,n\rangle.
$$

{\rk{Step 3}} Add a single stable letter
 $\sigma$ that conjugates the
free subgroup of $E_2$ generated by the $s_i$ to the
free subgroup of $E_2$ generated by the $\tau_i$:
$$
E_3 = \langle E_2, \sigma\mid \sigma^{-1}s_i\sigma
=\tau_i, \ i=0,\dots,n\rangle.
$$
At this stage  we have a group in which all of the
generators except $\sigma$ are conjugate to $g_0$. In
particular, every finite quotient  of $E_3$
is cyclic.

{\rk {Step 4}} Because no power of $a_0$ 
lies in either of the subgroups of $E_2$ generated
by the $s_i$ or the $\tau_i$, the normal form theorem
for HNN extensions implies that $\{a_0,\sigma\}$
freely generates a free subgroup of $E_3$.

We define $\widehat G$ to be an amalgamated free product of
 two copies of $E_3$,
$$
\widehat G = E_3\ast_F \-{E_3},
$$
where $F=F(x,y)$ is a  free group of rank
two; the inclusion into $E_3$ is  
$x\mapsto a_0$ and $ y\mapsto \sigma$, and the
inclusion into $\-E_3$ is  
$x\mapsto \-\sigma$ and $ y\mapsto \-a_0$.
All of the generators of 
 $\widehat G$ are conjugate to a power of
either $g_0$ or $\-{g_0}$, and therefore cannot survive
in any finite quotient. In other words, $\widehat G$ has
no finite quotients.
\qed


The following lemma enables us to gauge the geometry
of the embeddings in the preceding construction.

\proclaim{2.2 Lemma} Let $G$ be a 
group with finite 
generating set $\Cal A$, where no $a\in\Cal A$
represents $1\in G$.  
\roster
\item In any HNN extension of $G$  with finitely
many stable letters $s_0,\dots, s_n$, the free subgroup
generated by $S=\{s_0,\dots, s_n\}$ is isometrically
embedded with respect to $\Cal A\cup S$. 
If $\langle a \rangle \subset G$ is isometrically embedded
and has trivial intersection with the amalgamated subgroups
of $s_i$ then $\gp\{a,s_i\}$ is isometrically embedded in
the HNN extension.
\item If $H\subset G$ is isometrically embedded
with respect to $\Cal A$,  then
$H$ is also isometrically embedded in any
isometric amalgamation involving $G$ as a vertex
group (provided the amalgamation is isometric
with respect to the same generating set $\Cal A$).
\item Let $g\in G\smallsetminus\{1\}$. The cyclic subgroups of 
$G\ast \langle t \rangle$ generated
by $t$, by $[g,t]$,  and by each 
$(at)$ with $a\in\Cal A$,
are all isometrically
embedded  with respect to the choice of generators
 $\Cal A^*=\{at, [g,t],t \mid a\in\Cal A\}$.
\endroster
\endproclaim

\demo{Proof} (1)  and (2)  follow  from the
normal form theorem  for graphs of 
groups \cite{Se}. 

The normal form theorem for free products tells
us that if we write $[g,t]^n$ as a word in the
generators $\Cal A\cup\{t\}$, then that word must contain at
least $2n$ occurences of  $t^{\pm 1}$.
Each of the  elements of $\Cal A^*$
contains at most two occurences of $t^{\pm 1}$, therefore
$d_{\Cal A^*}(1, [g,t]^n)=n$.

If a word over $\Cal A\cup\{t\}$  equals $(at)^n$ in $G\ast\<t\>$,
then its exponent sum in $t$  
must be $n$. Therefore, since each of the generators in $\Cal A^*$
has $t$-exponent sum $1$ or $0$, we have $d_{\Cal A^*}(1, (at)^n)=n$.
\qed

\demo{\bf 2.3 The Proof of Theorem A} 
We   follow the proof of (2.1).
What we must ensure is that at each stage the embedding 
which we described   can be performed by
means of an {\it isometric} amalgamation. 

First we choose a finite generating set $\Cal A$ for
$G_0 = G\ast T(n)$ so that $G\hookrightarrow G_0$
is an isometric embedding, and we fix an element $g\in G_0$
whose image is trivial in every finite quotient of $G_0$. Then
as generators for $G_1= G_0\ast\langle t \rangle$
we take $\Cal A^*:= \{at, [g,t],t \mid a\in\Cal A\}$.
Note the   difference with (2.1) --- we have included
$[g,t]$. Define $g_0=[g,t]$.

Lemma 2.2(3) assures us that the
amalgamations carried out in Step 1 of the proof
of (2.1) are along
isometrically embedded subgroups  provided that
we take all $p_i=1$. And parts (1) and (2) of  Lemma 2.2
imply that the amalgamations
carried out in  Steps 2, 3 and 4 of (2.1) are
also along isometrically embedded subgroups.
Thus we obtain the desired group $\widehat G\in \Cal G$ that
has no finite quotients.

We have the inclusions $G\subset G_0\subset G_1\subset \widehat G$.
The third inclusion was constructed to be an isometric embedding.
The first and second
inclusions are obviously isometric embeddings
with respect to natural choices of generators.
But it does not follow that $G\hookrightarrow \widehat G$
is an isometric embedding, because at the end of
 Step 0 of the   proof we switched
from the obvious set of generators for $G_1$ to
a less natural set that was suited to our purpose.
On the other hand,  for any finitely generated group $H$, 
the identity map between the metric spaces
obtained by endowing $H$ with
different word metrics is  bi-Lipschitz.
 Thus,  $G\subset\widehat G_0$
is a quasi-isometric embedding (with respect to any choice
of word metrics).
\qed

For future reference we note:

\proclaim{2.4 Lemma} The cyclic subgroups generated by
all of the stable letters introduced in the above
construction are isometrically embedded in $\widehat G$.
\endproclaim

\heading{The non-positively curved case}\endheading

The proof that we shall
give of  Theorem C is entirely self-contained
except that we do not prove the basic facts about non-positively
curved spaces that are listed (3.2). 
One could shorten the proof of Theorem C considerably by using the 
complexes constructed in
\cite{W3} or \cite{BM} in place of   Lemmas 3.3 and 3.5. However
those constructions are rather complicated, so we feel that
there is  benefit in presenting a more direct account.

The example given in (4.3(2))   shows 
that the class of groups which act properly and
cocompactly on spaces of non-positive curvature
does not satisfy the conditions of
Theorem A. Nevertheless, with appropriate attention
to detail, one can use  the blueprint of our proof of
 Theorem A to
prove Theorem C, and this is what we shall do.
First   we need to know that there exists a
compact non-positively curved 2-complex whose
fundamental group is not residually finite.  

\demo{\bf 3.1 Wise's Examples \cite{W1}}
Let
$$
T(n)= \langle a, b, t_a, t_b \mid [a,b]=1,\,
t_a^{-1}at_a= (ab)^n,\, t_b^{-1}bt_b=(ab)^n \rangle.
$$ 
In Section 1 we saw that if $n\ge 2$ then this group is not Hopfian
and therefore not residually finite.
$T(n)$ is the fundamental group of the non-positively curved
2-complex $X(n)$ that one constructs as follows:  
take  the  (skew) torus 
formed by identifying opposite sides  of a rhombus  with
sides of length  $n$ and small diagonal of length 1; 
the loops formed by the images of the sides of the rhombus
are labelled $a$ and $b$ respectively; to this torus  
attach   two tubes $S\times [0,1]$, where $S$ is a
circle of length $n$; one end of the first tube is attached to the
loop labelled $a$ and one end of the second tube
is attached to the loop labelled $b$; in each case
the other end of the tube   wraps $n$ times around
the image of the small diagonal of the rhombus.  

Any complex obtained by attaching tubes along local
geodesics in the above manner  is non-positively curved 
in the natural length metric (see 
[\Cite{BH}, II.11]). We shall  need
the following additional facts  
concerning metric spaces of non-positive curvature; see
\cite {BH} for details.

\proclaim{3.2 Proposition} Let $X$ be a compact, connected, geodesic
space of non-positive curvature. Fix $x\in X$.
\roster 
\item Each homotopy class in
 $\pi_1(X,x)$ contains   a unique shortest
loop based at $x$. This based loop is the unique local geodesic
 in the given  homotopy class.
\item  Each conjugacy class in
 $\pi_1(X,x)$ is represented by a closed geodesic in $X$
(ie~a locally
isometric embedding of a circle). In other words, every
loop in $X$ is freely homotopic to a closed geodesic
(which need not   pass through $x$). If two closed geodesics
are freely homotopic then they have the same length.
\item $\pi_1(X,x)$ is torsion-free.
\item Metric graphs are non-positively curved.
\item The induced path metric on the 1-point union of
two non-positively curved spaces is again non-positively
curved.
\item If $X$ is a compact non-positively curved space, $Z$
is a compact length space
and  $i_1, i_2 \co  Z\to X$ are locally isometric embeddings,
then, when endowed with the induced path metric,
the quotient of 
$X \cup (Z\times [0,L])$ by the equivalence relation generated
by $i_1(z) \sim (z,0)$ and $i_2(z)\sim (z,L)$ is non-positively
curved. Moreover, if $L$ is greater than the diameter of
$X$, then $X$ is isometrically embedded in the
quotient.
\endroster
\endproclaim

A particular case of (6) that we shall need
is where $X$ is the disjoint union of
spaces $X_1$ and $X_2$,  and
$Z$ is a circle. In this case
the quotient  is obtained by joining $X_1$ to $X_2$ with a
cylinder  whose ends are attached along  closed geodesics.

\proclaim{3.3 Lemma} There exists a compact,
connected,
non-positively curved 2-complex $K$ with basepoint $x_0\in K$
such that:
\roster
\item there is an element $g_0\in \pi_1(K,x_0)$ whose
image in every finite quotient of $\pi_1(K,x_0)$ is trivial;
\item $\pi_1(K,x_0)$ is generated by a finite set of elements
each of which is
represented by a closed geodesic  that passes through $x_0$
and has integer length;
\item $g_0$ is represented by a closed geodesic of length 1 
 that passes through $x_0$.
\endroster
\endproclaim

\demo{Proof} Let $X$ be a compact, connected, 2-complex
of non-positive curvature and let $g_0\in \pi_1X$
be a non-trivial element whose image in every finite
quotient of $\pi_1X$ is trivial (the spaces $X(n)$
of (3.1) give such examples).
We choose a point $x_0$   on a closed
geodesic that represents the conjugacy class of $g_0$.
Suppose that $\pi_1(X,x_0)$ is generated
by $\{b_1,\dots,b_n\}$, let $\beta_i$ be the
shortest loop based at $x_0$ in the homotopy
class $b_i$, and let   $l_i$ be the length
of $\beta_i$.  Let $l_0$ be the
length of the closed geodesic representing $g_0$. Replacing
$g_0$ by a proper power if necessary, we may assume that
$l_0>l_i$ for $i=1,\dots,n$.

Consider the following metric graph $\Lambda$:
 there are $(n+1)$ vertices
$\{v_0,\dots,v_n\}$ and $2n$ edges $\{e_1,\e_1,\dots,e_n,\e_n\}$;
the edge $e_i$ connects $v_0$ to $v_i$ and has length $(l_0-l_i)/2$;
the edge $\e_i$ is a loop of length $l_0$ based at $v_i$.
We obtain the desired complex $K$ by gluing $\Lambda$
 to $X$, identifying $v_0$ with $x_0$,
 and then scaling the metric by a factor of
$l_0$  so that the closed geodesic representing $g_0\in \pi_1(K,x_0)$
has length 1.

Let $\g_i\in\pi_1(K,x_0)$ be the element given by the geodesic
$c_i$ that traverses
$e_i$, crosses $\e_i$, and then returns along $e_i$, that
is $c_i=e_i\e_i\-e_i$, where the overline denotes reversed
orientation.
 Note that $\pi_1(K,x_0)$ is the free product of
 $\pi_1(X,x_0)$ and the free group generated by 
$\{\g_1,\dots,\g_n\}$. As generating set for
$\pi_1(K,x_0)$ we choose $\{b_i\g_i, b_i\g_i^2 \mid i=1,\dots,n\}$.

According to parts (4) and (5) of the preceding
proposition, $K$ has non-positive curvature.
Moreover, the concatenation of
any non-trivial locally geodesic loop in $X$, based at $x_0$,
and any non-trivial locally geodesic loop in $\Lambda$ based at $v_0$
is a closed geodesic in $K$. Thus $\b_ic_i$ and
$\b_ie_i\varepsilon_i^2
\-e_i$ are closed geodesics in $K$; the former has
length $2$ and the latter has length $3$;  the former
represents $b_i\g_i$ and the latter represents $b_i\g_i^2$.
\qed

\demo{\bf 3.4 The proof of Theorem C} 
Given a compact, connected,
non-positively curved space $X$ we must isometrically embed
it in a  compact, connected,
non-positively curved space $\overline X$ whose fundamental
group has no non-trivial finite quotients. Moreover the
embedding must be such that if $X$ is a complex
of dimension at most $n\ge 2$ then so is $\overline X$.
We give two constructions, the first in outline and the
second in detail.

{\bf First Proof}\stdspace  We form the
1-point union of $X$ with one of the complexes $X(n)$
described in (3.1) thus ensuring that
some element $g_0$ of the fundamental group has
trivial image in every finite quotient. We then
apply the construction of (3.3), gluing a metric
graph to our space to obtain a space $X'$ whose
fundamental group is generated by elements
represented by closed geodesics that pass through
a basepoint on a closed geodesic representing $g_0$. 
To complete the
proof   one follows the argument of Lemma 3.5
  with $X'$ in
place of $K$ (taking the cylinders attached to
be sufficiently long so that $X$ is isometrically
embedded in the resulting space,  3.2(6)).

{\bf Second Proof}\stdspace Choose a finite set of generators
for  $\pi_1X$, and let $c_1,\dots,c_N$
be   closed geodesics in $X$ representing
the conjugacy classes of these elements. 
Lemma 3.5 gives
a compact non-positively
curved 2-complex $K_4$ whose fundamental group has no
finite quotients; fix  a closed geodesic $c_0$ in $K_4$.
Take $N$ copies of $K_4$  
 and scale the metric on the $i$-th copy so that the length
of $c_0$ in the scaled metric is equal to the length $l(c_i)$ of
$c_i$. Then glue the $N$ copies of $K_4$ to $X$ using cylinders
 $S_i\times [0,L]$ where $S_i$ is a 
circle of length $l(c_i)$; the ends of $S_i\times [0,L]$ 
 are attached by arc length
parametrizations of $c_0$ and $c_i$ respectively.
Call the resulting space $\overline X$.

Part (6) of  (3.2) assures us that $\overline X$ is
non-positively curved, and if the length $L$ of the gluing tubes
is sufficiently large then the natural embedding
$X\hookrightarrow \overline X$ will be an isometry.
\enddemo

It remains to construct $K_4$.

\proclaim{3.5 Lemma} There exists a compact non-positively
curved 2-complex $K_4$ whose fundamental group has no
finite quotients.
\endproclaim

\demo{Proof} 
Let $K$ be as in (3.3).
We mimic the argument of (2.1), with $\pi_1(K,x_0)$
in the r\^ole of $G_1$. At each stage we
shall state what the fundamental group of the
complex being constructed is; in each case this
is a simple application of the Seifert-van Kampen
theorem.

 Let $c_0$ be the
closed geodesic of length 1 representing $g_0$.
 Let $\{a_0,\dots,a_n\}$  
be the generators given by  3.3(2),   let
$\alpha_i$ be the closed geodesic through $x_0$ that represents
$a_i$, and
suppose that  $\alpha_i$ 
has length $p_i$.
 For each $i$, we glue to $K$ a cylinder
$S_{p_i}\times [0,1]$, where $S_{p_i}$ is a circle 
of length $p_i$, with 
basepoint $v_i$; one end of the cylinder is attached 
to $\alpha_i$ while the other end  
wraps $p_i$-times around 
$c_0$, and 
$v_i\times\{0,1\}$ is attached to $x_0$.
Let $K_1$ be the resulting complex. By the
Seifert-van Kampen theorem, $\pi_1(K_1,x_0)=E_1$,
in the notation of (2.1). Part (6) of (3.2)
implies that $K_1$
is non-positively curved.

The images in $K_1$ of the paths
$v_i\times [0,1]$ give  an isometric embedding into $K_1$ 
of the metric graph $Y$ that
has one vertex and  $n$ edges of length 1;
call the corresponding free subgroup $F_1\subset E_1$ 
(it is the subgroup generated by the $s_i$ in (2.1)).

Step 2 of (2.1) is achieved by attaching $n$   
cylinders of unit circumference 
$S_1\times [0,1]$ to $K_1$, the ends of
the $i$-th cylinder being attached to   $c_0$
and to the image of $v_i\times [0,1]$. The
resulting complex $K_2$ has $\pi_1(K_2,x_0)=E_2$.
As in the previous step,
 the free subgroup $F_2\subset E_2$   generated
by the basic  loops that run along the new cylinders
is the $\pi_1$-image of an isometric embedding    
$Y\to K_2$. (This $F_2$ is the subgroup generated by the
$\tau_i$ in (2.1).)

To achieve Step 3 of (2.1), we now glue $Y\times [0,L]$
  to $K_2$ by attaching the ends according to the
isometric embeddings that realize the embeddings
$F_1,F_2\subset  \pi_1(K_2,x_0)$. This gives us
a compact non-positively curved complex $K_3$ with
fundamental group $E_3$ (in the notation of
(2.1)). Let $v$ be the vertex of $Y$, observe that
$v\times\{0,L\}$ is attached to $x_0\in K_3$, and let
 $\sigma\in\pi_1(K_3,x_0)$ be the homotopy class of
the loop $[0,L]\to K_3$ given by
$t\mapsto (v,t)$.

We left open the choice of $L$, the length of the mapping
cylinder in Step 3, we now specify that it should
be $p_0$, the length of the geodesic representing
the generator $a_0$.
An important point to observe is that the angle at
$x_0$ between the image of $v\times [0,L]$ and any
path in $K_1\subset K_3$ is $\pi$. Thus
the free subgroup $\gp\{a_0,\sigma\}$ is the $\pi_1$-image
in $\pi_1(K_3,x_0)$ of an isometry from  the metric graph 
$Z$ with
one vertex (sent to $x_0$) and two edges of
length $L=p_0$. In fact, we have two such isometries  $Z\to K_3$,
corresponding to the free choice we 
have of which  edge of $Z$ to send to the image of $v\times [0,L]$.
We use these two maps to
 realize Step 4 of the construction on
(2.1): we apply part (6) of (3.2)
with $X$ equal to the disjoint union of
two copies of $K_3$  and with the two maps $Z\to 
K_3$ employed as the local isometries  $i_1, i_2$, the image of
one of the maps being  in each component of $X$.
The resulting space is the desired complex $K_4$.
\qed

By gluing non-positively curved orbi-spaces (in the
sense of Haefliger \cite{Ha}), or by performing
equivariant gluing,
 one can extend Theorem C to include
groups with torsion. We refer the reader 
to  [\Cite{BH}, II.11] for the technical
tools that make  this adaptation straightforward.

\proclaim{3.6 Theorem} If a group $G$ acts properly and
cocompactly by isometries on a \c space $Y$ then one can
embed $G$ in a group $\widehat G$ that acts  properly and
cocompactly by isometries on a \c space $\overline Y$ 
and has no proper subgroups of finite index. If $Y$ is
a polyhedral complex of dimension $n\ge 2$ then so is $\overline Y$.
\endproclaim

Since the group $G$   need not be torsion-free,
(3.6) shows in particular that there exist
compact non-positively
curved orbihedra,   with finite local
groups, that are not finitely covered by any
 polyhedron (where `covered' refers to
covering in the sense of orbispaces and `polyhedron'
means an orbihedron whose local groups are trivial).
We close our discussion of non-positively curved spaces
with an explicit example to illustrate this point.
The first examples
of this type were discovered by my student Wise \cite{W2}, and
the following example is essentially contained in his work. 
   
\example{\bf 3.7 A semihyperbolic group that is
not virtually torison-free}

In  the hyperbolic plane    
$\Bbb H^2$ we consider 
 a regular quadrilateral $Q$ with vertex angles $\pi/4$.
 Let $\alpha$ and $\beta$ be 
hyperbolic translations that identify the opposite sides
of $Q$. Then $Q$ is a fundamental domain for the action
of   $G=\gp\{\alpha, \beta\}$; the commutator 
 $[\alpha,\beta]$
 acts as a rotation through  $\pi$ at one vertex of $Q$, and
away from the orbit of this vertex the action of $G$ is
free.
Thus the quotient orbifold $V=\Bbb H^2/G$ is a torus
with one singular point, and at that
singular point the local group is $\Bbb Z_2$.

Let $X(n)$ and $T(n)$ be as in (3.1) and fix a 
closed geodesic
$c$
in the homotopy class of a non-trivial element $g_0$ in the kernel
 of a self-surjection $T(n)\onto T(n)$.  We scale the metric
on $X(n)$ so that this geodesic has length $l=|\alpha |=|\beta |$. Then we
take a copy of $X(n)$ and consider the orbispace $\overline V$
 obtained
by gluing it to $V$ using a tube $S_l\times [0,1]$  one
end of which is glued to $c$ and the other end of
which is glued to the image in $V$ of the axis of  $\alpha$. 

$\overline V$
inherits the structure as a (non-positively
curved) orbihedron in which the
only singular point is the original one; at this
singular point
the local structure is as it was in $V$. The
 fundamental group $\widehat G$ 
of  
$\overline V$ is $G\ast_\Bbb Z T(n)$, where the amalgamation
identifies $g_0\in T(n)$ with $\alpha\in G$. Now,
$g_0$ has trivial image in every finite quotient of $T(n)$,
therefore $[\alpha, \beta] = [g_0, \beta]$  has trivial image in 
every finite quotient of $\widehat G$.  
It follows that
$[\alpha, \beta]$, which has order two, lies in
every   subgroup of $\widehat G$ that has finite index.

In the case $n=2$,  the  
group $\widehat G$ has the following presentation:
$$
\langle a,b,s,t,\a,\beta\mid \a=[s^{-1}
(ab) s, b],\, [a,b]=[\a,\beta]^2=1,\, 
t^{-1}bt=s^{-1}as=(ab)^2\rangle.
$$
\endexample

\heading{Isoperimetric inequalities}\endheading

Isoperimetric inequalities for finitely presented groups
$G=\<\Cal A\mid\Cal R\>$
measure the complexity of the word problem. If a word  
$w$ in the free group $F(\Cal A)$
represents the identity in   $G$, then there
is an equality
$$ w=\prod_{i=1}^N x_i^{-1}r_ix_i$$
in $F(\Cal A)$, where $r_i\in \Cal R^{\pm 1}$. 
Isoperimetric inequalities give upper bounds on the
integer $N$ in a minimal such expression. 
The bounds are given as a function of the length of $w$,
 and the function $f_G\co \Bbb N\to\Bbb N$
giving the optimal bound is called
the {\it Dehn function} of the presentation.
If there is a constant $K>0$ such that the
functions $g,h\co \Bbb N\to \Bbb N$ satisfy
$g(n) \le K\,  h(Kn) + Kn $, then one writes $g\preceq h$.
It is not difficult to show (see \cite{Al} for     
example) that the Dehn functions
of different   finite presentations 
of a fixed group  are $\simeq$
equivalent, where $f \simeq g$ means that $f\preceq g$
and $g\preceq f$.

As an alternative measure of complexity for
the word problem, instead of trying to bound
the integer $N$ in the above equality one
might seek to bound  
the length of the conjugating elements $x_i$. In this case the
function giving the optimal bound is called the
{\it isodiametric function} of the group, which
we write $\Phi_G(n)$.  Again, this
function
is $\simeq$ independent of the chosen presentation 
(see \cite{Ge}). 

We refer the reader to \cite{Ge} for more information and
references concerning
Dehn functions
and isodiametric functions and their (useful) interpretation
in terms of the geometry of van Kampen diagrams.

\proclaim{4.1 Proposition} If $G$ is   an isometric
amalgam  of a finite collection 
$\{G_i \mid i\in I\}$ of finitely
presented groups, then 
the Dehn function $f_G(n)$ of $G$ is $\preceq  n^2 + n\,
\max_i f_{G_i}(n) $.
\endproclaim
\demo{Proof} A diagrammatic version of the 
proof is given in (4.3(3)), here we present a
more algebraic proof.

 By definition, $G$ is
the fundamental group of a finite graph of
groups. For the sake of notational convenience we shall assume that
there are no loops in the graph of groups under
consideration. The proof in
the general case is entirely
similar but notationally cumbersome.

Thus we have a finite tree with vertex set $I$ and a set
of edges $\Cal E\subset I\times I$. At the vertex indexed $i$
the  vertex group is 
 $G_i$. Let $H_{i,j}$ be the edge group associated to $(i,j)\in\Cal E$.
By definition, (1.3),
 there are finite generating
sets $\A_i$ for the
$G_i$  and subsets $\Cal B_{i,j}\subset \A_i$  with specified
bijections $\phi_{i,j}\co \Cal B_{i,j} \rightarrow \Cal B_{j,i}$
for each $(i,j)\in \Cal E$; the set
  $\Cal B_{i,j}$ generates $H_{i,j}$, each of
the inclusions $ H_{i,j}\hookrightarrow G_i$ is isometric
with respect to these choices of generators, and 
 $\phi_{i,j}=  \phi_{j,i}^{-1}$.

We fix finite
 presentations $\<\Cal A_i\mid\Cal R_i\>$  for the $G_i$. Then,
$$G \cong \<\Cal A\mid
\Cal R, \  \phi_{i,j}(b)=b,\, \forall b\in
 \Cal B_{i,j} \>,$$ where $\Cal A =\coprod_i \Cal A_i,
\, \Cal R = \coprod_i \Cal R_i$, and
 $(i,j)$ runs over $\Cal E$

Let $W$ be a word in the generators $\Cal A$. Suppose that $W$
is identically equal to  a product
$u_1\dots u_m$, where each $u_k$ is a   word 
over one of the alphabets
$\A_{i(k)}$ and each $\A_{i(k)}\neq  \A_{i(k+1)}$. Under
these circumstances $W$ is said to have {\it alternating
length $m$}. The normal form theorem for amalgamated
free products \cite{LS}
(or more generally  graph products  \cite{Se}) ensures that
this notion of length is well-defined. It also tells us
 that if $W=1$ in $G$
then  at least one   of the subwords  $u_k$
is equal in $G_{i(k)}$ to a word $\omega$
in the generators $\Cal B_{i(k),i(k\pm 1)}$.
Because $H_{i(k),i(k\pm 1)}$ is isometrically
embedded in $G_{i(k)}$, we can replace $u_k$ by
$\omega$ without increasing the length of $W$. This
can be done at the cost of applying at most 
$f_{G_{i(k)}}(2|u_k|)$
relations. We apply $|\omega|$ relations
to replace each letter $b$ of $\omega$ 
with $\phi_{i(k), i(k\pm 1)}(b)$.
Then, without applying any
more relations, we group $\omega$ together
with the neighbouring word $u_{k\pm 1}$. The
 net effect of this operation is to
reduce the alternating length of $W$ without
increasing its actual length. By repeating this
operation fewer than
$|W|$ times we can replace $W$ by a word $W'$ with
$|W'|\le |W|$ that involves letters from only  one
of the alphabets $\Cal A_i$. Since $W'$ represents the identity in
$G_i$, we can then reduce $W'$
to the empty word by applying at most $f_{G_i}(|W'|)$ relators
from $\Cal R_i$. 

The total number of relators applied in the reduction of
$W$ to $W'$ is fewer than $m |W| + m\, \max_i f_{G_i}(|W|)$, where
$m$ is the alternating length of $W$.
Therefore the total number of relators that we had to apply
in reducing $W$ to the empty word was less than
$|W|^2 + |W|\, \max_i f_{G_i}(|W|)$.
\qed

\proclaim{4.2 Corollary} The class of groups that
satisfy a polynomial isoperimetric inequality is closed
under the formation of isometric amalgamations along finitely
generated subgroups.
\endproclaim

\remark{4.3 Remarks} 

(1)\stdspace If instead of considering isometric amalgamations
we considered the fundamental groups of graphs of
groups in which the edge groups were only
quasi-isometrically embedded, then the above
proof would break down at the point where we
noted that 
$|W'|\le |W|$. In fact Proposition 4.1 would be false
under this weaker hypothesis: consider the
Baumslag-Solitar groups for example.

(2)\stdspace Let $D$ be the direct product of the free group
on $\{a,b\}$ and the free group on  $\{c,d\}$. Let $L=\gp\{ac, bc\}$.
For a suitable choice of generators,
 $L$ is isometrically embedded in $D$. It is shown
in \cite{B} and \cite{BH} that  $D\ast_LD$ has a cubic Dehn
function, whereas $D$ has a quadratic Dehn function.
Thus, in general, isometric amalgamations may increase the
polynomial degree of  Dehn functions.

(3)\stdspace
The proof of (4.1) can be recast as an induction
argument in which one proves that
the area of a minimal van Kampen diagram for $W$
is\nl  $m(\max_i f_{G_i}(|W|) +|W|)$, where $m$ is the
alternating length of $W$. This admits a simple 
geometric proof which we shall now sketch.

Draw a circle labelled by $W$, divide  it into
$m$ subarcs according to the decomposition of
$W$ as an alternating word. Maintaining the notation
established in the proof of (4.1), we
draw a chord in the disc connecting the endpoints of  
the circular arc  labelled by  $u_k$.
We label the chord by a geodesic word $\omega
\in\Cal B_{i(k), i(k\pm 1)}^*$ that is equal to $u_k$
in $G$. We  fill
the subdisc with boundary labelled $u_k\omega^{-1}$  
using a minimal-area
van Kampen diagram over the given presentation of
$G_{i(k)}$. We then  attach to the chord
 labelled $\omega$  faces corresponding to relators
of the type $\phi_{i(k), i(k\pm 1)}(b)$; the effect of this
is
to replace
$\omega$ by the corresponding
word in the generators $\Cal B_{i(k\pm 1),i(k)}$.
By induction, we may fill the
remaining   
subdisc  with a van Kampen diagram of area no greater than
$(m-1)(\max_i f_{G_i}(|W|)+|W|)$. We may choose $u_k$
so that $2|u_k| \le |W|$, and hence $|u_k|+|\omega|\le |W|$.
Therefore the area of the
whole diagram is no greater than
 $m\, (\max_i f_{G_i}(|W|)+|W|)$, completing the induction. 
\endremark

A simple induction on alternating length,
in the manner of (4.3(3)), allows one to show that 
(with respect to the finite presentations considered
in (4.1)) every
null-homotopic word $W$
of alternating length $m$ bounds a van Kampen diagram
in which every vertex can be joined to the basepoint of
the diagram by a path in the 1-skeleton 
that has length at most $|W| + \max_i\Phi_{G_i}(|W|)$.
Thus:

\proclaim{4.4 Proposition} If $G$ is an isometric
amalgam  of a finite collection 
$\{G_i \mid i\in I\}$ of finitely
presented groups, then  the isodiametric function $\Phi_G(n)$ of 
$G$ is $\preceq\max_i \Phi_{G_i}(n)$.
\endproclaim

\demo{\bf 4.5 The Proof of Theorem B} 
Given an infinite finitely presented
group $G$, we replace it by $G\ast \Bbb Z$.
This does not change the Dehn function or the
isodiametric function of $G$  but it allows us to assume that
$G$ is generated by a finite set of
elements $\{a_i,\dots,a_r\}$ such that each $\<a_i\>$ is
isometrically embedded in $G$ (see 2.2(3)).

The fundamental group $S$ of any of the spaces $\overline X$ yielded
by Theorem C will satisfy a quadratic isoperimetric
inequality and a linear isodiametric
inequality [\Cite{BH}, III]. At the level of $\pi_1$, the proof of
Theorem C was exactly parallel to that of (2.1), so Lemma 2.4
implies that   $S$ contains an isometrically embedded
infinite cyclic subgroup $\<s\>$. 

The group $\widehat G$ whose existence is asserted in
Theorem B is obtained by taking
an amalgamated free product of $G$
and $m$ copies
of $S$:  the cyclic subgroup $\<s\>$ in the
$i$-th copy of $S$ is identified with $\<a_i\>\subset G$.
In other words, $\widehat G$ is the fundamental group
of a tree of groups in which there is one vertex of valence $m$,
with vertex group $G$, and $m$ vertices of
valence 1, each with  vertex group   $S$; each edge 
group is infinite cyclic and the generator of the $i$-th
edge group is mapped to $s\in S$ and $a_i\in G$.

Proposition 4.1 tells us that the Dehn function of $\widehat G$
is $\preceq nf_G(n)$, and Proposition 4.4 tells us that
the isodiametric function of $\widehat G$
is no worse than that of $G$.
\qed

\Addresses
\recd

\enddocument